\documentclass[a4paper,12pt]{article}
\usepackage{amsmath}
\usepackage{graphicx}
\usepackage{rotating}
\usepackage{wasysym}
\DeclareGraphicsExtensions{.pdf,.png,.jpg}
\usepackage{chicago}
\usepackage{klett}
\usepackage{def_the}

\renewcommand{\theequation}{\mbox{\arabic{section}.\arabic{equation}}}
\renewcommand{\thefigure}{\arabic{section}.\arabic{figure}}
\renewcommand{\thetable}{\arabic{section}.\arabic{table}}
\renewcommand{\footnoterule}{\rule{14.8cm}{0.3mm}\vspace{+1.0mm}}
\renewcommand{\baselinestretch}{1.0}

\title{}
\author{Eckhard Platen}
\begin{document}
\thispagestyle{empty} \vspace*{1.0cm}

\begin{center}
{\LARGE\bf Entropy-Maximizing  Dynamics  \\
	\vspace*{0.5cm}
	of  Continuous Markets}
\end{center}

\vspace*{.5cm}
\begin{center}

{\large \renewcommand{\thefootnote}{\arabic{footnote}} {\bf Eckhard
Platen}\footnote{University of Technology Sydney,
  School of Mathematical and Physical Sciences, and \\Finance Discipline Group}$^{,}$
\footnote{Australian National University, Canberra, College of Business and Economics}}
\vspace*{2.5cm}

\today

\end{center}

\begin{minipage}[t]{13cm}
By assuming the existence of the growth optimal portfolio (GOP), the stationarity of GOP-volatilities, and the maximization of  relative entropy, the paper applies the benchmark approach to the modeling of the  long-term dynamics of   continuous markets.  It reveals  conservation laws, where the GOP is shown to follow a time-transformed squared Bessel process of dimension four.  Moreover, it predicts the convergence of the averages of the   GOP-volatilities with respect to the driving independent Brownian motions toward a common level. 

\end{minipage}
\vspace*{0.5cm}

{\em JEL Classification:\/} G10, G11

\vspace*{0.5cm}
{\em Mathematics Subject Classification:\/} 62P05, 60G35, 62P20
\vspace*{0.5cm}\\
\noindent{\em Key words and phrases:\/}   growth optimal portfolio, market price of risk, entropy maximization, conservation law,   
  squared Bessel process, Radon-Nikodym density.  

\vspace*{0.5cm}
\noindent{\em Acknowledgements:\/} The author would like to express his gratitude for
\noindent receiving valuable suggestions on the  paper by Martin Schweizer, Mark Craddock, Kevin Fergusson, Martino Grasselli, Juri Hinz, Hardy Hulley, Erik Schl\"ogl, Michael Schmutz,  and Stefan Tappe.

	\newpage
	\section{Introduction}\label{section.intro}

	A recent mathematically rigorous account of what  one may call classical mathematical finance theory  is given in \citeN{Jarrow22}.  The   {\em No Free Lunch with Vanishing Risk} (NFLVR) no-arbitrage condition of \citeN{DelbaenSc98},  also considered in \citeN{DelbaenSc94}  and \citeN{DelbaenSc06}, provides the theoretical basis for risk-neutral pricing and represents the core assumption of the classical mathematical finance theory; see, e.g., \citeN{Jarrow22}. As cited and summarized in \citeN{Jarrow22}, enormous progress has been made in the understanding  of financial markets. However, every theory has  limits to its applicability. \\
	 As a first  example that indicates a  limitation of the classical mathematical finance theory, one can point out its denomination dependence. When modeling by assuming NFLVR with the  
	   savings account of the domestic currency as num\'eraire, there can still exist
	    free lunches with vanishing risk from the perspective of a foreign market participant who prices risk-neutrally  by employing the foreign savings account  as num\'eraire. This  is unsatisfactory not only from a theoretical point of view. It has also the practical consequence that risk-neutral pricing  works  economically only for short-term contingent claims. \\
	     A second example, where classical mathematical finance theory has been unsatisfactory, is given by the asymptotic singularity of the risk-neutral pricing measure for typical models when the time to maturity increases to infinity. This fact hinders  realistic long-term modeling, pricing, and hedging. \\
	    A third example emerges from the fact that 
	    the  widely followed investment strategy for retirement suggests  investing when young in a diversified portfolio of risky securities and  gliding over time into   more and more investments in  fixed-income securities. This strategy  has been successfully followed by    generations. It takes  advantage  of the  higher average growth rate of  diversified portfolios of risky securities compared to the average growth rate of the savings account.   Classical mathematical finance theory suggests 
	        buying at any time during the accumulation phase the number of  units of the savings account that one wishes to own  at the retirement date  and  holding these units  until  retirement. This  investment strategy is theoretically not consistent with the above-mentioned successfully practiced  financial planning strategy.  \\
	  The above three examples  show that  an extension of the classical mathematical finance theory is needed if one aims to overcome the mentioned limitations.
	The current paper derives new results in a direction that provides solutions for the above-described shortcomings  by applying the author's benchmark approach; see \citeN{Platen06ba} and \citeN{PlatenHe06}.  The    {\em growth optimal portfolio} (GOP) plays a central role in the benchmark approach, where it is called {\em benchmark} and any security is said to be {\em benchmarked} when  denominated in units of the GOP.  The GOP is   interchangeably  called the Kelly portfolio, expected logarithmic utility-maximizing portfolio, or num\'eraire portfolio; 
	see  \citeN{Kelly56}, \citeN{Merton71}, \citeN{Long90}, \citeN{Becherer01}, \citeN{Platen06ba}, \citeN{KaratzasKa07}, and \citeN{HulleySc10}. Due to the properties of the logarithmic function, the  weights of the GOP are denomination-independent, which forms the basis for the  denomination independence of the benchmark approach.\\ 
	 	As with any theory,	  the proposed application of the benchmark approach must depend on simplifying assumptions that may be not quite true  but permit the successful  engineering of solutions to typical practical tasks in an idealized market setting. These assumptions should only be based on empirial evidence and mathematical reasoning. \\The  paper considers  a {\em basis market}, which  represents an idealized, nonredundant, continuous  market with independent benchmarked {\em basis security accounts} as primary security accounts. The independent benchmarked basis security accounts  can be obtained as self-financing portfolios of the benchmarked primary security accounts of a  general, nonredundant, continuous     market, as will be shown in Appendix D. A basis market  has its  stochastic independent   benchmarked  basis security accounts and the constant benchmarked GOP as primary security accounts. 
	 	  \\
	 	   In \citeN{PlatenRe20} it has been shown that a well-diversified total return stock index is a good proxy for the GOP of the stock market. A wide range of the econometric literature identified the stylized empirical fact that  volatility processes of stock indexes exhibit equilibrium properties; see, e.g., \citeN{Engle82}. This literature observed also  the leverage effect, where a stock index and its volatility are negatively correlated; see, e.g., \citeN{Black76b}.
	 Based on these stylized empirical facts, the paper models the volatilities of the GOP in the denominations of  basis security accounts as  scalar, stationary diffusions, which have a stationary density  and  random initial values  distributed according to their stationary density; see, e.g., Section 4.5 in \citeN{PlatenHe06}. \\With such a stationary  volatility process the average of the instantaneous growth rate of a basis security account-denominated GOP becomes  a constant. Consequently, a basis security account-denominated GOP can be normalized by the exponential  of a respective {\em activity time} when the    derivative of the activity time equals  the  average instantaneous growth rate of the  basis security account-denominated GOP. 
	 	To model parsimoniously the volatility of the resulting {\em normalized GOP} as a scalar  diffusion, the simplest model appears to be  a local volatility function model,   pioneered, e.g., by \citeN{BreedenLi78}, \citeN{DermanKa94a}, and \citeN{Dupire94}. Therefore, the paper models this volatility  as a  function of the normalized GOP. \\
	 	To identify the local volatility function, the current paper maximizes the  entropy of the stationary density of the normalized GOP. The concept of entropy was  pioneered by \citeN{Boltzmann1866}  and has been crucial for the mathematical identification of realistic dynamics  in many areas.  For a stochastic dynamical system, its relative entropy,  see, e.g., \citeN{Kullback59}, is  typically increasing 
	 	and moves over time toward its maximum.
	 	 The candidate pricing density for  basis security accounts as num\'eraire is characterized by the  benchmarked basis security accounts.  Its relative entropy  with respect to the stationary density of the normalized GOPs  appears to represent the relative entropy  of a basis market. In summary, the paper makes the following three  key assumptions:\\
	 	
	 \noindent	{\bf{A1:}}
	  For a basis  market  the GOP exists, and its  benchmarked basis security accounts are independent scalar diffusions. \\

	 	\noindent {\bf{A2:}} 
	 	 A normalized  GOP   forms a 
	 	 stationary scalar diffusion and its volatility is a function of its value. \\

	 	\noindent {\bf{A3:}} The market  
	 	maximizes  the relative entropy of  the   candidate pricing  density  with respect to the stationary joint  density  of the  normalized  GOPs. \\

	 	    The first assumption A1 is  about the existence of the GOP and    represents an  intuitive and easily verifiable {\em no-arbitrage condition}   because \citeN{KaratzasKa07} and \citeN{KaratzasKa21} have shown that  the existence of the GOP is equivalent to their {\em No Unbounded Profit with Bounded Risk} (NUPBR)  condition. This no-arbitrage condition is weaker than the NFLVR  condition of \citeN{DelbaenSc98}. Under the NFLVR condition, the GOP always exists. By relying on a weaker assumption,  the benchmark approach accesses  a  richer   modeling world than the classical mathematical finance theory covers. \\

	  The maximization of the relative entropy is known to be equivalent to the  minimization of the information rate; see \citeN{Kullback59}. 
	 	Consequently, the resulting  entropy-maximizing market dynamics  does not leave any room for exploitable information. 
	 	In this sense, the third key assumption A3 postulates some kind of  {\em market efficiency}, which is different from most other  notions of market efficiency that have been  discussed in the  literature; see, e.g., \citeN{Fama70}. \\
	 	
	 	Conservation laws  guide and simplify the dynamics of complex dynamical systems in many areas. According to \citeN{Noether18}, the maximization of a Lagrangian  in the presence of Lie-group symmetries for the underlying model dynamics leads to the identification of 
	 	conservation laws; see   \citeN{Olver93} and  \citeN{KosmannSc18}.   The entropy-maximizing stationary dynamics of a normalized GOP  turns out to have   Lie-group symmetries and respective conservation laws emerge. 
	 	 The  entropy-maximizing  dynamics of normalized GOPs  are those of time-transformed square root processes, see \citeN{RevuzYo99}.  Their conserved dimension is four and their conserved logarithmic mean is zero, which fixes their arithmetic mean  at a respective  level. The entropy-maximizing dynamics  for the GOP in basis security account denomination exhibits a levergage effect and has Student-t distributed log-returns  with four degrees of freedom. Both properties coincide with  stylized empirical facts.\\ 
	 	 
	 	Benchmarked basis security accounts turn out to be the inverses of  time-transformed squared Bessel processes of dimension four, which are known to be  strict supermartingales and not  true martingales; see, e.g., \citeN{RevuzYo99}. Consequently, for  a given  basis security account as num\'eraire, the candidate pricing density turns out to be a strict supermartingale. 
	 	 Classical mathematical finance theory postulates the pricing density of a risk-neutral pricing measure to form a true martingale. However, under the entropy-maximizing dynamics this assumption is not satisfied when the savings account is a basis security account. 
	 	 Therefore, when the savings account is a basis security account, the  widely applied risk-neutral pricing rule, a cornerstone of classical mathematical finance theory, 
	 	  	 has under the above three key assumptions  no theoretical basis as the only appropriate pricing method and more economical prices and hedges  exist; see, e.g., \citeN{Platen06ba}.\\

	 	  As a result of  the relative entropy maximization it turns out that the time derivatives of the activity times become equal to each other. This makes  the market prices of risk, which are the volatilities of the basis security account-denominated GOPs, on average equal. This fact can be interpreted as the outcome of  the  search by buyers and sellers  for  appropriate prices  in their decision making. \\
	 	  
	  To exploit the power and elegance of continuous stochastic calculus,  the current paper  focuses on continuous financial markets with continuous trading,  continuous security paths,  instantaneous investing and borrowing,   short sales with full use of proceeds, and   infinitely divisible securities. 
	    Since most of the resulting market properties are not very sensitive to  minor market frictions,  the current paper   assumes    no transaction costs. \\
	    
	     The paper is organized as follows:  Section 2  introduces  a  continuous basis market 
	     with  flexible local volatility functions  for the volatilities of normalized GOPs.
	      Section 3  specifies the local volatility functions by maximizing the entropy of the stationary densities of the normalized GOPs. 
	       Section 4 maximizes the relative entropy of the  candidate pricing density. 
	       The first three appendixes prove respective theorems. The fourth appendix shows that the predicted dynamics of basis markets determine also the  dynamics of more general markets.	
	\section{Basis Market}
\subsection*{Basis Security Accounts}

	 Since for general markets, one cannot expect a set of simple relationships to emerge, we introduce  a {\em basis market} that allows us to derive under the three key assumptions   its entropy-maximizing dynamics. We explain in Appendix D, by transforming a basis market into a more general market,  that these dynamics determine the dynamics of more general markets.  \\ 
	 The  modeling throughout the paper is performed on a filtered probability space $(\Omega,\mathcal{F},\underline{\cal{F}},P)$, satisfying the usual conditions; see, e.g., \citeN{KaratzasSh88} and \citeN{KaratzasSh98}.  We introduce  the vector process ${\bf{W}}=\{ {\bf{W}}_t =(W^1_t,\dots,W^{n}_t)^\top,t \in [0,\infty)\}$ of the $n\in\{1,2,...\}$   independent  standard {\em Brownian motions}  driving the dynamics of the 
	 securities in a market. We make the important assumption that the number of these Brownian motions is {\em irreducible} for the market dynamics we model, which means that there is no way to capture the fluctuations of the securities  with a lower number than $n$   independent Brownian motions. The filtration $\underline{\cal{F}}$ $=(\mathcal{F}_t)_{t \in [0,\infty)}$ models the evolution of information relevant to the model.
	The information available at time $t \in [0,\infty)$ is captured by the sigma-algebra $\mathcal{F}_t$. 
	For matrices and vectors $ \bf{x}$ we 
	denote by $ {\bf{x}}^\top$ their transpose.  Moreover, $\b1=(1, \dots, 1)^\top$ is a vector, and we write $\bf 0$ for a zero matrix or vector, where the dimensions follow from the context.  \\

By the first key assumption A1 the GOP of the market is assumed to exist. 
 As can be seen 
 in  equation \eqref{bmportfolio} in  Appendix D, benchmarked self-financing portfolios in continuous markets are local martingales.     For $j\in\{1,...,n\}$, the $j$-th {\em basis volatility} is defined as the negative volatility of 
the $j$-th {\em benchmarked basis security account}  $\hat B^j_t$, which is the
   local martingale
that	satisfies the driftless SDE

\begin{equation}\label{hatBj}
\frac{ d \hat B^{j}_{t}}{\hat B^{j}_{t}}=
-\theta^j_tdW^j_t
\end{equation}
for $t\in[0,\infty)$ with initial value $\hat B^{j}_{0}=1$. It is only driven by the $j$-th Brownian motion $W^j$ and its volatility is  by the second key assumption A2 assumed to be a stationary scalar diffusion, which we assume for our manipulations to be square-integrable and predictable. According to the second key assumption A2, the $j$-th basis volatility is only impacted by randomness generated by the $j$-th Brownian motion $W^j_t$.    When holding  units of the $j$-th benchmarked basis security account, its owner is reinvesting in this account all  interest payments, dividends,  or  payments  that she or he     receives  or is obliged to pay. 
 For  the $j$-th {\em basis security account-denominated GOP} $(\hat B^j_t)^{-1}$, which is the inverse  of the $j$-th benchmarked basis security account,  one obtains by  application of the  It\^{o} formula   the SDE
\begin{equation}
\frac{ d(\hat B^j_t)^{-1}}{(\hat B^j_t)^{-1}}=
\theta^j_t(\theta^j_tdt+dW^j_t)
\end{equation}
for $t\in[0,\infty)$ for $j\in\{1,...,n\}$. 
 To form a nonredundant continuous market with  basis security accounts as primary security accounts, we add to the above-introduced $n$ stochastic benchmarked basis security accounts the constant benchmarked GOP $\hat B^{n+1}_t=1$ as $(n+1)$-th benchmarked basis security account. 
 
 \subsection*{Normalized GOP}

  For  $j\in\{1,...,n\}$, the logarithm of the $j$-th basis security account-denominated GOP $(\hat B^j_t)^{-1}$ satisfies by  application of the  It\^{o} formula  the SDE
  \begin{equation}\label{lnhatBj}
   d(\ln(\hat B^j_t)^{-1})=
  \theta^j_t(\frac{\theta^j_t}{2}dt+dW^j_t)
  \end{equation}
  for $t\in[0,\infty)$. 
   By the second key assumption A2,  we define for $j\in\{1,...,n\}$ and $t\in[0,\infty)$ the $j$-th {\em normalized  GOP} 
   \begin{equation}\label{Y}
   Y^j_{\tau^j_t}=\frac{(\hat B^j_t)^{-1}}{e^{\tau^j_t}},
   \end{equation} as the normalized  $j$-th basis security account-denominated GOP.  
   We assume it to represent a stationary scalar diffusion process in the    $j$-th {\em activity time}
    \begin{equation}\label{tau}
    \tau^j_t=\tau^j_0 +
   \ a^jt
    \end{equation}
    for $t\in[0,\infty)$, where   $\tau^j_0$ denotes the  {\em j-th  initial activity time}.. We denote by $q^j_\infty$ the stationary density of the $j$-th normalized GOP $Y^j_.=\{Y^j_{\tau}, \tau\in[\tau_0,\infty)\}$  and by ${\bf{E}}^{q^j_\infty}\left(.\right)$ the respective expectation. We  call the latter  the respective {\em average}  to indicate that it represents the expectation with respect to the stationary density $q^j_\infty$.  The   constant $j$-th {\em activity} \begin{equation}a^j=\frac{d\tau^j_t}{dt}={\bf{E}}^{q^j_\infty}\left(\frac{(\theta^j_.)^2}{2}\right)>0\end{equation} is the derivative of the $j$-th activity time and equals the average of the stationary instantaneous growth rate $\frac{(\theta^j_t)^2}{2}$ of the   $j$-th basis security account-denominated GOP.   
  Note, since $Y^j_{\tau^j_.}$ is assumed to be a stationary process,  its initial value
  \begin{equation}\label{Y0}
  Y^j_{\tau^j_0}=e^{-\tau^j_0}
  \end{equation}
       is distributed according to the stationary density  $q^j_\infty$ of the $j$-th normalized GOP.
     \subsection*{Basis Volatility}
      The $j$-th   basis volatility   $\theta^j_t$ is the volatility of the $j$-th basis security account-denominated GOP $(\hat B^j_t)^{-1}$, and is also the volatility of the  $j$-th normalized GOP $Y^j_{\tau^j_t}$. By the second key assumption A2, the $j$-th  basis volatility  is assumed to be a stationary process.  Therefore, it cannot be  a function of the  GOP because the latter is not a stationary process; see  SDE \eqref{e.4.1} in Appendix D.  To keep the model parsimonious, 
      the second key assumption A2 postulates  that the $j$-th basis volatility $\theta^j_t$ evolves in the $j$-th activity time  as the value of a local volatility function of the $j$-th normalized  GOP $Y^j_{\tau^j_t}$.   Since the $j$-th normalized GOP is assumed to evolve in the $j$-th activity time,  
      and by avoiding   a potential division by zero in the derivation of the resulting formula for the stationary density (see \eqref{qjy}) of the normalized GOP,  a suitable   ansatz for the 
   $j$-th basis volatility  is given by the formula
\BE\label{M}
\theta^{j}_{t}=\sqrt{ \frac{a^j}{ \phi^j(Y^j_{\tau^j_t})}}.
\EE
   Here $\phi^j(.)$ denotes a  flexible, at zero infinitely often continuously differentiable function of the $j$-th normalized GOP value that  is assumed to allow the  manipulations  we  perform.  
 With this  flexible function  it follows by the application of \eqref{Y},  \eqref{M}, \eqref{hatBj}, and the It\^{o} formula that the  $j$-th normalized GOP satisfies in the $j$-th activity time the SDE
\BE \label{dY0}
dY^j_{\tau^j}=Y^j_{\tau^j}\left(\frac{1}{\phi^j(Y^j_{\tau^j})}-1\right)d\tau^j +Y^j_{\tau^j}\sqrt{\frac{1}{\phi^j(Y^j_{\tau^j})}}d\bar W^j_{\tau^j}
\EE
for $\tau^j\in[\tau^j_0, \infty)$ with the initial value $Y^j_{\tau^j_0}$. Here the stochastic differential $d\bar W^j_{\tau^j_t}=\sqrt{a^j}d W^j_{t}$ characterizes the  Brownian motion $\bar W^j_.$ that evolves in the $j$-th activity time $\tau^j_t$ for $t\in[0,\infty)$. \\ It is well-known that the stationary density $q^j_\infty$ of the $j$-th normalized GOP in its activity time is the  solution of the following stationary Fokker-Planck equation (see, e.g., Chapter 4 in \citeN{PlatenHe06}):
\begin{equation}
\frac{d \left(q^j_\infty(y)y(\frac{1}{\phi^j(y)}-1)\right) }{dy}-\frac{1}{2}\frac{d^2 \left(q^j_\infty(y)y^2\frac{1}{\phi^j(y)}\right) }{(dy)^2}=0,
\end{equation} which is a  second-order ordinary differential equation. Its solution is given by the formula 
\BE \label{qjy}
q^{j}_\infty(y)=\frac{C_j\phi^j(y)}{y^2 }\exp \left\{2\int_{0}^{y}\frac{u(1-\phi^j(u))}{u^2}du\right\}
\EE
for $ y\in(0,\infty)$ and some constant $C_j>0$. 
     The  average  of the logarithm of the normalized GOP is a core quantity in our modeling and we parametrize it 
     as the  logarithmic average   
     \BE \label{stalogmean}
     {\bf{E}}^{q^j_\infty}(\ln(Y^j_{.}))=\zeta_j \EE
  by using the flexible parameter $\zeta_j \in (-\infty, \infty)$.  The average  of the $j$-th normalized GOP  $Y^j_.$ represents the  reference level of the stationary process $Y^j_.$ and  determines the average value of the GOP at a given time. We parametrize it as the arithmetic average
     \BE\label{stamean}
     {\bf{E}}^{q^j_\infty}(Y^j_{.})=\bar Y_j
     \EE 
    by using a flexible parameter $\bar Y_j\in[0,\infty)$ for $j\in\{1,...,n\}$. \\ 
    Since the function $\phi^j(.)$ is   flexible, the above-resulting     volatility model covers a wide range of local volatility function models, including those   studied in \citeN{DermanKa94a}, \citeN{Dupire94},  \citeN{BreedenLi78}, and    
     \citeN{Cox96}. 
      A  deeper theoretical reason for a particular choice of any of these local volatility function models has been so far missing in the literature. The current paper will provide a theoretical basis for a realistic local volatility function model. \\
    In reality,  feedback effects in trading activity make the activity processes  appear as fast-moving stationary  processes, as observed, e.g., in \citeN{PlatenRe19} and \citeN{GuyonLe22}. In long-term modeling the  stochasticity of a fast moving stationary activity process  generates only  minor second-order effects in pricing and hedging. Therefore, the  randomness of activity processes  
     is, for simplicity, neglected in the current paper. However, it can be explained and modeled  in the way as shown in \citeN{PlatenRe19}. \\
     The following notion of a basis market will  allow us to derive simple relationships for its dynamics:   
           \begin{definition}
            	A {\em  basis market} is a  continuous market with  $n\in\{1,2,...\}$ driving independent Brownian motions, which is formed by the above-introduced $n+1$  independent benchmarked basis security account processes $\hat B^1,...,\hat B^{n+1}$. 
      \end{definition} 
  \section{Entropy Maximization}
 \subsection*{Radon-Nikodym Density}     
   When  the $j$-th basis security account,  $j\in\{1,...,n\}$, is the denominating  security, the inverse $(\hat B^j_t)^{-1}$ of the $j$-th benchmarked basis security account, see \eqref{hatBj}, is the GOP in this denomination. The  respective candidate  pricing measure  is  characterized by
    the {\em j-th Radon-Nikodym density} 
\begin{equation}\label{lambda}
   \Lambda_{(\hat B^j)^{-1}}(t)=  \hat B^j_t=  \frac{1}{e^{\tau^j_t}Y^j_{\tau^j_t}}, 
   \end{equation}
  which equals the $j$-th benchmarked basis security account 
  for     $t\in[0,\infty)$, see \eqref{hatBj} and \eqref{Y}. 
   The above Radon-Nikodym density may not be a probability density because under the first key assumption A1 and by \eqref{bmportfolio} in Appendix D the Radon-Nikodym density process $\Lambda_{(\hat B^j)^{-1}}$ is only a local martingale and {\em not} required   to form a true martingale. 
    \subsection*{Relative Entropy}
   Under the assumptions A2 and A3  the function $\phi^j(.)$ is specified when the stationary density $q^j_\infty$ is found.
     To specify the stationary density, the third key assumption A3 postulates that for  $t\in(0,\infty)$ the market maximizes the relative entropy of its Radon-Nikodym density  $\Lambda_t=\prod_{j=1}^{n}\Lambda_{(\hat B^j)^{-1}}(t) $  with respect to the stationary joint density $q_\infty=\prod_{j=1}^{n}q^j_\infty$.
   The following common definition of the {\em relative entropy} can be found,  e.g., in \citeN{Kullback59}:
    \begin{definition} For a basis market  and  $t\in(0,\infty)$, the {\em relative entropy} ${\cal{H}}(q_\infty,\Lambda_t)$ of the Radon-Nikodym density  $\Lambda_t$ with respect to the stationary joint density $q_\infty$ 
    	 is defined as
    	\BE\label{Hpq}
    	{\cal{H}}(q_\infty,\Lambda_t)=- 
    	\int q_\infty (y)\ln(\Lambda_t)dy= {\cal{H}}(q_\infty) -{\cal{D}}(q_\infty,\Lambda_t),
    	\EE
    	where \BE \label{entropy} {\cal{H}}(q_\infty)=	{\cal{H}}(q_\infty,q_\infty)=-
    	\int q_\infty (y)\ln(q_\infty(y))dy\EE
    	denotes the {\em entropy} of $q_\infty$, and
    	$
    	{\cal{D}}(q_\infty,\Lambda_t)
    	$
    	the {\em  Kullback-Leibler divergence} of $q_\infty$ from $\Lambda_t$.  
    \end{definition}
    The Kullback-Leibler divergence is  defined  through \eqref{Hpq} and \eqref{entropy}.  It is assumed  that all  quantities in the above definition are assumed to be   well defined componentwise for the manipulations we  perform, where the  integrations are running from $0$ to $\infty$. 
    \subsection*{Entropy-Maximizing  Dynamics}
  The third key assumption A3 postulates the maximization of the relative entropy ${\cal{H}}(  q_\infty,\Lambda_t)$. Due to the independence of the normalized GOPs and the  additive structure of the entropy for the involved joint densities, this maximization can be  performed in several steps:
  First, we 
  maximize for $j\in\{1,...,n\}$ the  entropy ${\cal{H}}(q^j_\infty)$ of the stationary probability density $q^j_\infty$ of $Y^j_{.}$ under 
  the constraints \eqref{stalogmean} and \eqref{stamean}. This maximization identifies the function $\phi^j(.)$, and therefore, for the j-th normalized GOP $Y^j_.$ in the $j$-th activity time $\tau^j_t$  its entropy-maximizing  stationary density $ \bar q^j_\infty$. 
   The proof of the following result is given in  Appendix A:  
  \begin{theorem}[Normalized GOP Dynamics]\label{YJT}
   For a basis market   with maximized entropy ${\cal{H}}(\bar q^j_\infty)$ 
    the function $\phi^j(.)$ takes the form
    \BE
   \phi^j(y)=\frac{y}{ \bar Y},
   \EE
   the respective stationary density is the gamma density
   \begin{equation}
   \bar q^{j}_\infty(y)=\left(\frac{2}{\bar Y}\right)^2y \exp\{-\frac{2y}{\bar Y}\}
   \end{equation}
   with conserved four degrees of freedom,
   	 conserved logarithmic average
 	 \begin{equation}\label{conserved}
 	 {\bf{E}}^{\bar q^{j}_\infty}(\ln(Y^j_{.}))=0,
 	 \end{equation} and conserved  arithmetic average
 	 \begin{equation}\label{barY'}
 	 \bar Y=2e^{\gamma_E-1},\end{equation}
 and the $j$-th normalized GOP satisfies the SDE
 \begin{equation}\label{dY2}
 dY^j_{\tau^j_t}=(\bar Y-Y^j_{\tau^j_t})d\tau^j_t+ \sqrt{  Y^j_{\tau^j_t}\bar Y} d\bar W^j_{\tau^j_t}
 \end{equation}	 for $t \in[0,\infty)$, with random initial value $Y^j_{\tau^j_0}$ distributed according to  $\bar q^{j}_\infty$ 
  	for $j\in\{1,...,n\}$.
  \end{theorem} 
  When the above entropy is  maximized, the   $j$-th normalized GOP  follows in the $j$-th activity time $\tau^j_t$
   a square root process   of dimension four; see, e.g., \citeN{RevuzYo99}. 
    The constant $\gamma_E$ in the equation \eqref{barY'} is the  {\em Euler-Mascheroni constant} with an approximate value of $\gamma_E \approx 0.5772$; see \citeN{Havil03}. The $j$-th basis security account denominated GOP  follows a time-transformed squared Bessel process of dimension four; see \citeN{RevuzYo99}.\\ 
       Since 
        the   entropy of the stationary joint  density of the normalized GOPs is employed as Lagrangian and this stationary density has Lie-group symmetries, as shown  in \citeN{CraddockPl04} and in Chapter 4 in \citeN{BaldeauxPl13}, Noether's Theorems, see \citeN{Noether18}, predict that conservation laws must  exist. 
         The  dimension four of the square root process is conserved, the logarithmic  average of the normalized GOPs is conserved at the  level zero, and the arithmetic average is conserved at the level $\bar Y$. These conservation laws arise mathematically  from the model structure,  and the first-order conditions of the  constrained entropy maximization in the proof of Theorem \ref{YJT}. 
 They simplify considerably the  model  and one has only to estimate as parameters the initial activity times and the activities to fit the model to historical data.\\
 The  resulting transition probability density $\bar q^j_t$ of the value of the $j$-th normalized GOP at time $t$ is a highly tractable non-central chi-square density, see equation (8.7.44) in \citeN{PlatenHe06}, which is the solution of the respective Fokker-Planck equation. It has  explicit formulas for many of its functionals (see  \citeN{BaldeauxPl13} for  an extensive collection  of such formulas).\\
  The obtained dynamics of squared   volatilities of benchmarked basis security accounts are those of the inverses of  time-transformed square root processes of dimension four. Such volatility dynamics are  known in the  literature as those of  3/2-volatility-type models. The derived 3/2-volatility  model was   first  suggested by \citeN{Platen97d} based on its links to the GOP.   The high tractability of  3/2  volatility-type models  was  observed in \citeN{Heston97}, who   suggested 3/2-volatility models as potential  option pricing models.\\ The entropy-maximizing dynamics of a basic security account-denominated GOP is captured by the {\em minimal market model}, which was  first suggested in \citeN{Platen97d} and named in   \citeN{Platen01a}; see also \citeN{HulleySc10}. The volatility  of  basis security account-denominated GOP dynamics is  negatively correlated to the GOP value, which matches the stylized empirical fact  known as the leverage effect; see, e.g., \citeN{Black76b}.\\ 
 Under the entropy-maximizing dynamics,  squared volatilities of GOPs in basis security account denomination  have  inverse gamma densities  with four degrees of freedom as  stationary densities. 
   Therefore,  when estimating the probability density of  log-returns of portfolios that approximate a GOP in the denomination of a basis security account 
    as if these were independent random variables,   a particular normal-mixture density must be expected. 
     The respective normal-mixture density is a {\em Student-t density} with four degrees of freedom. As previously mentioned, a well-diversified total return stock index   represents a good proxy for the GOP of the respective stock market. Several  independent empirical studies on log-return densities of stock indexes,  including \citeN{MarkowitzUs96b}, \citeN{HurstPl97d}, and \citeN{PlatenRe08e}, have  identified  the Student-t density with about four degrees of freedom as the probability density that matches with   high significance those of log-returns of well-diversified stock indexes.\\
     Benchmarked basis security accounts turn out to be the inverses of time transformed squared Bessel processes of dimension four, which are known to be strict supermartingales and not martingales; see \citeN{RevuzYo99}. Therefore, when the savings account is a basis security account, the risk-neutral pricing density is not a probability density.
      The real-world pricing formula of the benchmark approach, see Chapter 9 in \citeN{PlatenHe06} and \citeN{DuPl16}, provides the minimal possible price for a contingent claim. Major economic benefits have been identified in comparison to risk-neutral pricing when using the entropy-maximizing market dynamics in  real-world pricing and hedging of long-term contingent claims; see, e.g., \citeN{Fergusson20}, \citeN{SunZhPl21}, and \citeN{FergussonPl22}.   The  mentioned papers  cover also the case of  a contingent claim that pays one unit of the savings account at maturity. They answer  positively  the  question of whether  financial planning  strategies can be derived in a  rigorous manner that  are more economical than the strategies suggested by the classical mathematical finance theory.\\
                Due to the mentioned strict supermartingale property of benchmarked basis security accounts, the  entropy-maximizing market dynamics predict  the presence of \textquoteleft cheap snacks' and \textquoteleft free thrills' (see \citeN{LoewensteinWi00a}), \textquoteleft arbitrage amounts' (see \citeN{Platen02g}), \textquoteleft free lunches with vanishing risk' (see \citeN{DelbaenSc06}), \textquoteleft asset bubbles' (see  \citeN{JarrowPoSh07a} and \citeN{Jarrow22}), and \textquoteleft money market bubbles' (see \citeN{BaldeauxIgPl18}). \\ 
               \section{Activity Equilibrium}
      \subsection*{Market Activity}  
       Under the revealed dynamics of the normalized GOPs we have a gamma density as their stationary density with four degrees of freedom and mean $\bar Y$. Therefore,  it follows, see, e.g., equation (8.7.42) in \citeN{PlatenHe06}, that
        \begin{equation}
       {\bf{E}}^{\bar q^j_\infty}(\frac{1}{Y^j_{\tau^j_t}})=\frac{2}{\bar Y},
        \end{equation} 
        which yields the $j$-th activity at time $t$ as one half of the average of the $j$-th squared basis volatility, that is,
       \begin{equation}
       	a^j=\frac{1}{2}{\bf{E}}^{\bar q^j_\infty}(\frac{a^j\bar Y}{Y^j_{\tau^j_t}})=\frac{1}{2}{\bf{E}}^{\bar q^j_\infty}((\theta^j_t)^2)
       \end{equation} 
       for $j\in\{1,...,n\}$ and $t\in[0,\infty)$. So far, we interpreted the activities as exogenously given quantities. Let us relax this assumption by allowing the activities to become adjusted by the market. 
      More precisely, instead of the activities  (the averages of half the squared market prices of risk)  representing some exogenously given constants, we  request only  that their average, the  {\em  market activity}
         \BE\label{marketactivity}
         a=\left(\frac{1}{n}\sum_{j=1}^{n}\sqrt{a^j}\right)^2
         \EE
        is an exogenously given constant.  
         The relative entropy maximization can be interpreted as  the  search by the market for the   appropriate average level of  market prices of risk. 
         In Appendix B the following outcome of this search is derived:
         \begin{theorem} [Activity Theorem]\label{AJ}
         	For a  basis market, where the  entropy\\ $ {\cal{H}}(\bar q_\infty)$ is maximized and only the market activity $a$ is fixed but not the individual activities, the maximum relative entropy \begin{equation} {\cal{H}}(\bar q_\infty,\Lambda_t)
         	=-nat 
         	 \end{equation} 
         	 emerges  when  all individual activities  equal  the market activity.
         \end{theorem}
        The above result shows that the relative entropy maximization  drives  the activities toward the same level. 
          We conceptualize this by saying: Two benchmarked basis security accounts are said to be in  {\em activity equilibrium} if both have the same activities. When all $n$ stochastic  benchmarked basis security accounts   are  in activity equilibrium with each other,  the  basis  market  is said to be in activity equilibrium.  When  two initially separate basis   markets, where each is in its own activity equilibrium, become connected so that  the same market participants can invest in both markets, then the activities can be expected to start to move  spontaneously toward a new  market activity for the connected   market.
          One notes that relative entropy-maximization equalizes the average growth rates of the stochastic  benchmarked basis security accounts.\\ 
          
           When all activities converge to zero, a special state for a basis market  is approached.  In this case, it follows  from the limiting behavior of benchmarked basis security accounts that
          the benchmarked  basis security accounts $\hat B^1_t,...,\hat B^{n+1}_t$ become asymptotically constants.         
          \begin{definition}
          	We say that a basis market is in {\em activity equilibrium  } when   
          	the activities are all equal and nonnegative.
          \end{definition}
          
          \subsection*{Denomination in the Basis Portfolio}
         When denominating the securities in units of a basis security account, only the respective driving Brownian motion attracts a nonzero market price of risk in the SDEs for portfolios in this denomination.
         To construct a denominating portfolio where all driving Brownian motions attract a nonzero market price of risk, we form the  benchmarked {\em basis  portfolio} $\hat S^0_t$  satisfying the SDE
         \begin{equation}\label{dS0t}
         \frac{d \hat S^0_t}{\hat S^0_t}=\sum_{j=1}^{n}\frac{d \hat B^j_t}{\hat B^j_t}
         \end{equation} for $t\in[0,\infty)$ with $\hat S^0_t=1$. The basis portfolio goes long with equal weight $1$ in each of the stochastic benchmarked basis security accounts and invests   the remaining weight  $1-n$  in the benchmarked GOP. 
         The inverse of $\hat S^0_t$ represents the GOP denominated in units of the basis portfolio and satisfies the SDE
          \begin{equation}\label{aggrGOP}
          \frac{d (\hat S^0_t)^{-1}}{(\hat S^0_t)^{-1}}=\sum_{j=1}^{n} \theta^j_t( \theta^j_tdt+dW^j_t)
          \end{equation} 
         for $t\in[0,\infty)$. 
          As we have seen, the activities are  equal when the basis market is in activity equilibrium. Therefore, the denomination in the basis portfolio is the one where  a market in activity equilibrium has the same average for the market prices of  risk. The basis portfolio of a basis market is a special portfolio. For a basis market that is only driven by one Brownian motion, the basis portfolio equals the first basis security account.
          The GOP in the denomination of the basis portfolio aggregates in its volatility the market prices of risk with respect to all $n$ driving independent Brownian motions.
           Its dynamics in this denomination are    characterized by the following result, which is derived in Appendix C:
          \begin{theorem}\label{GPLA}
          	For a given   basis market  in activity equilibrium, 
          	the  {\em normalized    basis portfolio-denominated GOP }
          	 \begin{equation}
          	Y_{\tau_t}=\frac{(\hat S^0_t)^{-1}}{e^{\tau_t}},
          	\end{equation} 
                    satisfies	in {\em aggregate activity time}
          	\begin{equation}\label{bartau}
          	\tau_t= \tau_0 +
          	 n at,
          	\end{equation}
                    	 the  SDE
          	\begin{equation}\label{dY4}
          	dY_{\tau_t}=(\bar Y-Y_{\tau_t})d\tau_t+ \sqrt{  Y_{\tau_t}\bar Y} d\bar W_{\tau_t},
          	\end{equation}
          	with initial value $Y_{\tau_0}=e^{-\tau_0}$, where   	
          	\begin{equation}
          	d\bar W_{\tau_t}=\sqrt{\frac{Y_{\tau_t}}{ a\bar Y}}\sum_{k=1}^{n}\theta^k_t dW^k_t
          	\end{equation}
          denotes the stochastic differential of the {\em aggregate Brownian motion} $\bar W_{\tau_t}$ that is evolving in the aggregate activity time $\tau_t$	for $t\in[0,\infty)$. 
          \end{theorem}
          This theorem shows that not only does the GOP in the denomination of every basis security account follow the dynamics of a time-transformed squared Bessel process of dimension four. Also,  in the denomination of the  basis portfolio the  GOP has these dynamics.  Important is the observation  that the derivative of its aggregate activity time and, therefore, its average growth rate, equals $n$ times the market activity.  
  	\section{Conclusion}
  	 The  paper identifies the entropy-maximizing long-term dynamics of the growth optimal portfolio and the market prices of risk of  continuous financial markets.   The predicted dynamics   match well   empirical evidence and question the existence of an equivalent risk-neutral probability measure. 
  	  Extensions of the derived results toward markets with   jumps, transaction costs, discrete trading, and stochastic activity processes, as well as, applications to long-term pricing and hedging of contingent claims in real markets, represent further steps in the proposed research direction.\\
  
     \appendix
     \textwidth14.5cm \textheight9in
     \topmargin0pt
     \renewcommand{\theequation}{\mbox{A.\arabic{equation}}}
     \renewcommand{\thefigure}{A.\arabic{figure}}
     \renewcommand{\thetable}{A.\arabic{table}}
     \renewcommand{\footnoterule}{\rule{14.8cm}{0.3mm}\vspace{+1.0mm}}
     \renewcommand{\baselinestretch}{1.0}
     \pagestyle{plain}
     \section*{Appendix A: Proof of Theorem \ref{YJT}}\setcA \setcB
     Since the normalized GOPs in the denominations of respective  basis security accounts are independent, we can perform the maximization of the joint entropy in several steps.
     Under the  constraints \eqref{stalogmean} and \eqref{stamean} we maximize  first for $j\in\{1,...,n\}$ the entropy ${\cal{H}}( q^{j}_\infty)$ of the stationary probability density $q^{j}_\infty$ of the $j$th normalized GOP $Y^j_{\tau^j_.}$, according to the formulas \eqref{entropy} and \eqref{Hpq}, which means that we maximize the  Lagrangian
     \begin{equation*}
     {\cal{L}}(q^{j}_\infty)=- 
     \int_{0}^{\infty} q^j_\infty (y)\ln(q^j_\infty (y))dy+\lambda_0 \left(\int_{0}^{\infty}q^{j}_\infty(y)dy-1\right)\end{equation*}\BE +\lambda_1\left(\int_{0}^{\infty}y q^{j}_\infty(y)dy-\bar Y_j\right)
     +\lambda_2\left(\int_{0}^{\infty}\ln(y)q^{j}_\infty(y)dy-\zeta_j\right),
     \EE
     where $\lambda_0, \lambda_1, \lambda_2$ are Lagrange multipliers. $ {\cal{L}}(q^{j}_\infty)$ is maximized when its Fr\'echet derivative $\delta {\cal{L}}(q^{j}_\infty)$, i.e., the first variation of ${\cal{L}}(q^{j}_\infty)$ with respect to admissible variations of $q^{j}_\infty$, becomes zero. This implies
     \BE
     \delta {\cal{L}}(\bar q^{j}_\infty)=\int_{0}^{\infty}\left(-\ln(\bar q^{j}_\infty(y))+\lambda_0+\lambda_1 y+\lambda_2 \ln(y)\right)\delta \bar q^{j}_\infty(y) dy=0.
     \EE
     The solution of the above first-order condition is the gamma density 
     \BE \label{py}
     \bar q^{j}_\infty(y)=\exp \{\lambda_0+ \lambda_1 y+\lambda_2 \ln(y)\}
     \EE
    for $y\in (0,\infty)$ with the constraint \begin{equation}
     \int_{0}^{\infty}\exp \{\lambda_0+ \lambda_1 y+\lambda_2 \ln(y)\}dy=1,
     \end{equation}
    for the Lagrange multipliers $\lambda_0, \lambda_1,\lambda_2$. It
      has  $2(\lambda_2+1)$ degrees of freedom and as another constraint the mean \BE {\bf{E}}^{\bar q^{j}_\infty}(Y^j_{.})=\frac{\lambda_2+1}{-\lambda_1}=\bar Y_j.\EE \\On the other hand, the SDE for the $j$-th normalized index is given  by \eqref{dY0}.  
   Consequently, the stationary density $q^{j}_\infty(y)$ of the $j$-th normalized GOP satisfies the Fokker-Planck equation  with the drift and diffusion coefficient functions of the SDE \eqref{dY0}. This yields  the stationary density $q^{j}_\infty(y)$ in the form given in \eqref{qjy}. The latter must equal the above-identified gamma density. By setting equal both expressions for the stationary density of the normalized GOP, the  function $\phi^j(y)$ is no longer flexible, and we should find its mathematical description.
       The Weierstrass Approximation Theorem says that each continuous function can be approximated on a bounded interval by polynomials.  When using polynomials for characterizing $\phi^j(y)$ and searching for a match of the stationary density \eqref{qjy} with the gamma density in \eqref{py},  one finds that only the polynomial
     \BE\label{psi1}
     \phi^j(y)=\frac{y}{ \bar Y_j}
     \EE
     provides such a match. 
     This yields for the process $Y^j_.$ the 
     stationary density
     \BE\label{qY}
      q^{j}_\infty(y)=\left(\frac{2}{\bar Y_j}\right)^2y \exp\{-\frac{2y}{\bar Y_j}\}.
     \EE
     This is a gamma density with four degrees of freedom and mean $\bar Y_j$.
     Due to the constraint \eqref{stalogmean}, we have the entropy of the stationary density
     \BE\label{Hq'}
     {\cal{H}}( q^{j}_\infty)=-{\bf{E}}^{ q^{j}_\infty}(\ln(Y^j_{.}))=\ln(\frac{2}{\bar Y_j})+  \gamma_E-1=-\zeta_j,
     \EE
    which is the logarithmic average of the stationary density. Thus, the average level
    $
     \bar Y_j=2e^{\gamma_E-1+\zeta_j}
     $
     maximizes  the entropy $ {\cal{H}}( q^{j}_\infty)$. The latter becomes fully maximized for the value 
     \BE\label{zeta}
     \zeta_j=0,
     \EE
     yielding the maximum possible entropy value  \begin{equation}
     {\cal{H}}(\bar q^{j}_\infty)=-{\bf{E}}^{\bar q^{j}_\infty}(\ln(Y^j_{.}))=0.\end{equation} and the reference level \BE\label{barY}
     \bar Y_j=\bar Y=2e^{\gamma_E-1}
     \EE for $j\in\{1,...,n\}$.
      As a consequence of the identified stationary density, the originally flexible function $\phi^j(.)$ becomes specified.  This fixes the  dynamics of the $j$-th normalized index as that of a square root process of dimension four, satisfying the SDE \eqref{dY2}.   
      The resulting square root process is a stationary process when its initial value is distributed according to its stationary density. We note that its volatility  satisfies the originally imposed predictability and square integrability assumptions, 
       which proves Theorem \ref{YJT}.

      \textwidth14.5cm \textheight9in
           \topmargin0pt
      \renewcommand{\theequation}{\mbox{B.\arabic{equation}}}
      \renewcommand{\thefigure}{B.\arabic{figure}}
      \renewcommand{\thetable}{B.\arabic{table}}
      \renewcommand{\footnoterule}{\rule{14.8cm}{0.3mm}\vspace{+1.0mm}}
      \renewcommand{\baselinestretch}{1.0}
      \pagestyle{plain}

      \section*{Appendix B: Proof of Theorem \ref{AJ}}\setcA \setcB
     
      Using the entropy-maximizing stationary joint density $\bar q_\infty=\prod_{j=1}^{n}\bar q^j_\infty
      $ of the normalized GOPs, it remains  the maximization of    the relative entropy, which becomes the negative Kullback-Leibler divergence
      \BE
         {\cal{H}}( \bar q_\infty,\Lambda_t)=  -{\cal{D}}( \bar q_\infty,\Lambda_t)
      = \sum_{j=1}^{n}{\bf{E}}^{\bar q^j_\infty}(
       \ln(\frac{\hat B^j_t}{\hat B^j_0}))=-\frac{1}{2}\sum_{j=1}^{n}\int_{0}^{t}{\bf{E}}^{\bar q^j_\infty}((\theta^{j}_s)^2)ds \EE
      for  $t\in(0,\infty)$. 
        By \eqref{M}, \eqref{psi1}, \eqref{barY} and ${\bf{E}}^{\bar q_\infty}((Y^j_{.})^{-1})= \frac{2}{\bar Y^j}$, see, e.g., (8.7.42) in \citeN{PlatenHe06}, we obtain
      \begin{equation}\label{vk2}
      {\bf{E}}^{\bar q_\infty}((\theta^{j}_t)^2)=a^j\bar Y^j 
      {\bf{E}}^{\bar q_\infty}((Y^j_{.})^{-1})=2a^j.
      \end{equation} 
      Therefore, it follows 
      \BE {\cal{H}}( \bar q_\infty,\Lambda_t)=
      - 
      \sum_{j=1}^{n}a^jt. 
           \EE
      With the  market activity $a$ given in \eqref{marketactivity}
      we obtain
      \begin{equation*}
       {\cal{H}}( \bar q_\infty,\Lambda_t)
       =- \sum_{j=1}^{n}
      (\sqrt{a}+(\sqrt{a^j}-\sqrt{a}))^2t\end{equation*}
      \BE
      =- 
      n at 
      - \sum_{j=1}^{n}
      (\sqrt{a^j}-\sqrt{a})^2t.
      \EE
      Thus, the  relative entropy reaches  its maximum
          when
      \BE
      a^j=a
      \EE
      for all $j\in\{1,...,n\}$.
       This is the case when all activities   equal  the market activity, which proves Theorem \ref{AJ}. 
      
      \textwidth14.5cm \textheight9in
            \topmargin0pt
      \renewcommand{\theequation}{\mbox{C.\arabic{equation}}}
      \renewcommand{\thefigure}{C.\arabic{figure}}
      \renewcommand{\thetable}{C.\arabic{table}}
      \renewcommand{\footnoterule}{\rule{14.8cm}{0.3mm}\vspace{+1.0mm}}
      \renewcommand{\baselinestretch}{1.0}
      \pagestyle{plain}
      \section*{Appendix C: Proof of Theorem \ref{GPLA}}\setcA \setcB 
      
            The GOP $(\hat S^0_t)^{-1}$  in the denomination of the basis portfolio satisfies by \eqref{aggrGOP} the SDE
      \begin{equation}\label{barV*}
      \frac{d  (\hat S^0_t)^{-1}}{ (\hat S^0_t)^{-1}}=|\theta_t|\left( |\theta_t|d{ t}+d W_t\right),
      \end{equation}
      with {\em aggregate   GOP volatility}
      \begin{equation}\label{bartheta1}
      |\theta_t|=\sqrt{
      	\sum_{k=1}^{n}(\theta^k_t)^2},
      \end{equation}
     
      and  {\em aggregate  Brownian motion} $W_t$  
      with stochastic differential
      \begin{equation}\label{barW1}
      d W_t=\frac{1}{ |\theta_t|}\sum_{k=1}^{n}\theta^k_t dW^k_t
      \end{equation}
           for $t\in[0,\infty)$. 
      By parametrizing  the   squared aggregate GOP-volatility with respect to the  aggregate Brownian motion $W_t$ as \begin{equation}\label{theta2}
      |\theta_t|^2=\frac{\bar Yn a}{  Y_{\tau_t} },
      \end{equation}
      with   normalized  GOP 
      \begin{equation}
      Y_{\tau_t}=\frac{(\hat S^0_t)^{-1}}{e^{\tau_t}},
      \end{equation} 
          and aggregate activity time
      \begin{equation}
      \tau_t= \tau_0 +
       n at,
      \end{equation}
           one obtains for the normalized GOP the  SDE
      \begin{equation}\label{dY5}
      dY_{\tau_t}=(\bar Y-Y_{\tau_t})d\tau_t+ \sqrt{  Y_{\tau_t}\bar Y} d\bar W_{\tau_t},
      \end{equation}
      with $Y_{\tau_0}=e^{-\tau_0}$, where

      \begin{equation}
      d\bar W_{\tau_t}=\sqrt{na}d W_t
      \end{equation}
      for $t\in[0,\infty)$.  
       Consequently,  the normalized  GOP follows a square root process of dimension four in the aggregate activity time,
              which proves Theorem \ref{GPLA}.
      
       \textwidth14.5cm \textheight9in
             \topmargin0pt
       \renewcommand{\theequation}{\mbox{D.\arabic{equation}}}
       \renewcommand{\thefigure}{D.\arabic{figure}}
       \renewcommand{\thetable}{D.\arabic{table}}
       \renewcommand{\footnoterule}{\rule{14.8cm}{0.3mm}\vspace{+1.0mm}}
       \renewcommand{\baselinestretch}{1.0}
       \pagestyle{plain}
       \section*{Appendix D: Transformed Basis Market} \setcA \setcB 
       \subsection*{Market of Reference}
              Consider a basis market  with 
         $n\in\{1,2,...\}$ driving independent Brownian motions and $n+1$ benchmarked basis security accounts $\hat B^1_t,...,\hat B^{n+1}_t$. We can transform this basis market into a more general market by forming $n+1$  self-financing strictly positive benchmarked portfolios $\hat S^1_t,...,S^{n+1}_t$  and interpret these as the primary security accounts of the more general transformed market. For $j\in\{1,...,n+1\}$, the $j$-th primary security account is formed  by using the predictable and square-integrable weight process $\tilde \pi^j=(\tilde \pi^{j,1},...,\tilde \pi^{j,n+1})^\top$ when investing in the $n+1$ basis security accounts.   The  $j$-th benchmarked primary security account $\hat S^j_t$ satisfies the SDE
       \begin{equation} 
       \frac{d\hat S^j_t}{\hat S^j_t}=\sum_{k=1}^{n+1}\tilde \pi^{j,k}_t\frac{d\hat B^k_t}{\hat B^k_t}=-\sum_{k=1}^{n}\tilde \pi^{j,k}_t\theta^k_t dW^k_t
       \end{equation}
       for $t\in[0,\infty)$ with $\hat S^j_t>0$. 
       The volatility vector  at time $t$ of the $j$-th benchmarked primary security account $\hat S^j_t$ we denote by ${\bf{u}}^j_t=
       ( u^{j,1}_t,...,u^{j,n}_t)^\top$, where $u^{j,k}_t=-\tilde \pi^{j,k}_t\theta^k_t$ for all $j\in\{1,...,n+1\}$,  $k\in\{1,...,n\}$, and $t\in[0,\infty)$. Furthermore, we set   $u^{j,n+1}_t=1$ for $j\in\{1,...,n+1\}$. Let the   square matrix ${\bf{u}}_t=[u^{j,k}_t]_{j=1,k=1}^{n+1,n+1}$ combine the component volatilities    of the  benchmarked primary security accounts.   
       When the  matrix ${\bf{ u}}_t$ is invertible for all  $t\in[0,\infty)$,  the resulting market is called 
       a {\em  market of reference}. 
          It is straightforward to conclude that the above market of reference has the same GOP and the same market prices of risk as the underlying  basis market. 
       \subsection*{Growth Optimal Portfolio of a General Continuous Market}		    	
       Let us now study  a more general class of continuous markets which includes the above market of reference. We  consider a  market  with $n+1$ primary security accounts  that is denominated  in some strictly positive price process and driven 
       by the $n$  independent Brownian motions $W^1_t,...,W^n_t$. The  primary security accounts are given by the vector of strictly positive  values ${\bf{S}}_t=(S^1_t,..., S^{n+1}_t)^\top$ that are satisfying at time  $t \in [0,\infty)$ the $(n+1)$-dimensional  It\^{o}-vector SDE 
       \begin{equation} \label{e.2.1}
       \frac{d{\bf{S}}_t}{{\bf{S}}_t}=\mu_t dt +
       {\bf{\sigma}}^{}_t d{\bf{W}}_t
       \end{equation}
             with vector of positive initial values  ${\bf{S}}_0$. The dynamics of the primary security accounts are characterized by the  appreciation rate vector $\mu_t=(\mu^{1}_t,\dots,\mu^{n+1}_t)^\top$, and the volatility matrix   ${\bf{\sigma}}_t=[\sigma^{j,k}_t]_{j,k=1}^{n+1,n}$ for   $t \in [0,\infty)$. Here we write $\frac{d{\bf{S}}_t}{{\bf{S}}_t}$ for the $(n+1)$-vector of stochastic differentials $(\frac{dS^1_t}{S^1_t},...,\frac{dS^{n+1}_t}{S^{n+1}_t})^\top$.
                     The components of the  processes ${\bf \mu}
       $ and  
       $\sigma
       $ 
       are assumed to be 
       predictable and square-integrable. 
            A positive self-financing portfolio $S^\pi_t$  is described by its positive initial value $S^\pi_0>0$ and its  {\em weight process} $\pi=\{\pi_t=(\pi^1_t, \pi^2_t,..., \pi^{n+1}_t)^\top,t\in[0,\infty)\}$ quantifying the predictable and square-integrable fractions of wealth invested in the respective primary security accounts, where
       \begin{equation}
       \pi_t^\top {\bf 1}=1
       \end{equation}  for $t \in [0,\infty)$. The portfolio value process satisfies the SDE
       \begin{equation} \label{e.2.12}
       \frac{dS^\pi_t}{S^\pi_t}=\pi_t^\top\frac{d{\bf{S}}_t}{{\bf{S}}_t}=\pi_t^\top \mu_t dt +\pi_t^\top 
       {\bf{\sigma}}^{}_t d{\bf{W}}_t
       \end{equation}
             for $t \in [0,\infty)$. 
       A GOP in this market is a strictly positive self-financing portfolio $S^\pi_t$ with initial value $S^\pi_0=1$ that maximizes its instantaneous growth rate 
       \begin{equation}\label{gpi}
       g^\pi_t=\pi^\top_t  \mu_t-\frac{1}{2}\pi^\top_t  \sigma_t  \sigma^\top_t  \pi_t
       \end{equation} for $t \in [0,\infty)$. 
       This leads us to the  $(n+1)$-dimensional constrained quadratic optimization problem
       \begin{equation}\label{quop}
       \max \left\{g^\pi_t| \pi_t \in {\bf R}^{n+1}, \pi^\top_t  {\bf 1}=1\right\}
       \end{equation}
       for $t \in [0,\infty)$.	 To characterize its solution, we introduce the symmetric 
       matrix \BE  {\bf M}_t=\left( \begin{array}{cc} \sigma_t \sigma^\top_t & \b1  \\ \b1^\top & 0 \end{array} \right) 
       \EE
       and denote by $im\left(	{\bf M}_t\right) $ the image of this matrix; see, e.g., \citeN{GolubVa96}. 
       By Theorem 3.1 in \citeN{FilipovicPl09} one obtains  the following result:       
            	A GOP exists in the above market if and only if 
       	\begin{equation}\label{imM}
       	\left( \begin{array}{c} \mu_t \\ 1 \end{array} \right)\in im\left(	{\bf M}_t\right) 
       	\end{equation}
       	for all $t \in [0,\infty)$. In the case when the above condition is satisfied, the following statements emerge by the just mentioned theorem: \\ The value process of the GOP $S^{\pi^*}_t$ is unique and satisfies the SDE
       	\BE \label{e.4.1}
       	\frac{dS^{\pi^*}_t}{S^{\pi^*}_t}=\lambda^*_t dt +
       	{\bf v}^\top_t  ( {\bf v}_t dt+d{\bf W}_t) \EE
       	for $t \in [0,\infty)$, where its optimal weights  $\pi^*_t=(\pi^{*,1}_t,...,\pi^{*,n+1}_t)^\top$ together with its  {\em  risk-adjusted return} $\lambda^*_t$
       	represent  a solution to the  equation
       	\BE  \left( \begin{array}{cc} \sigma_t \sigma^\top_t & \b1  \\ \b1^\top & 0 \end{array} \right) \left( \begin{array}{c} \pi^*_t \\ \lambda^*_t 
       	\end{array} \right)=\left( \begin{array}{c} \mu_t \\ 1 \end{array} \right), \EE
       	       	with the vector of {\em GOP-volatility components}
              	\BE \label{e.4.2} {\bf v}_t=(v^1_t,...,v^n_t)^\top=\sigma^\top_t  \pi^*_t
       	\EE
       	       	and the   risk-adjusted return
       	\begin{equation}\label{lambda*}
       	\lambda^*_t=
       	(\pi^*_t)^\top\mu_t-{\bf v}_t^\top{\bf v}_t 
       	\end{equation}
              	for $t \in [0,\infty)$. Furthermore, the value of a strictly positive self-financing portfolio $S^\pi_t$ with weight process $\pi$ satisfies the SDE
       	\begin{equation}\label{portfolio}
       	\frac{dS^{\pi}_t}{S^{\pi}_t}=\lambda^*_t dt +\pi_t^\top
       	\sigma_t  ( {\bf v}_t dt+d{\bf W}_t), \end{equation}
       	and its benchmarked value $\hat S^\pi_t=\frac{S^\pi_t}{S^{\pi^*}_t}$  the SDE
       	\begin{equation}\label{bmportfolio}
       	\frac{d\hat S^{\pi}_t}{\hat S^{\pi}_t}={\bf u}^\pi_td{\bf W}_t \end{equation}
       	with volatility
       	\begin{equation}\
       	{\bf u}^\pi_t=\pi_t^\top
       	\sigma_t- {\bf v}_t  \end{equation}
       	for $t \in [0,\infty)$. 
            Under the condition \eqref{imM}, the  drift  coefficients of the SDEs for  benchmarked self-financing portfolios  in the above general market are always zero, which shows that benchmarked portfolios are local martingales. Nonnegative local martingales are supermartingales, which makes benchmarked nonnegative portfolios supermartingales.\\
      When the above market is a market of reference, one can assume that the relative entropy of the underlying basis market has been maximized. This gives the GOP-volatility components and market prices of risk in the market of reference the entropy-maximizing dynamics. 
       
\bibliographystyle{chicago}
\bibliography{my}

\newpage

 \end{document}